\documentclass[12pt,twoside]{amsart}\usepackage{amssymb}
\usepackage{enumerate}

\newtheorem{thmspec}{\relax}
\newtheorem{theorem}{Theorem}[section]
\newtheorem{thm}[theorem]{Theorem}
\newtheorem{lem}[theorem]{Lemma}

\newtheorem{prop}[theorem]{Proposition}
\newtheorem{defi}[theorem]{Definition}

\theoremstyle{definition}

\theoremstyle{remark}

\numberwithin{equation}{section}
\def \Bbb{\mathbb}

\def\onto{{\kern3pt\to\kern-8pt\to\kern3pt}}

\def\<{\langle}
\def\>{\rangle}
\def\|{{\ |\ }}

\def\onto{\twoheadrightarrow}

\def\-{\underline}

\def\inte{\operatorname{int}}

\def\R{\Bbb R}
\def\C{\Bbb C}

\def\X{\Bbb X}

\def\<{\langle}
\def\>{\rangle}

\catcode`\@=11
\def\serieslogo@{\relax}
\def\@setcopyright{\relax}
\catcode`\@=12

\title[Generalization of a theorem of Gonchar]
{Generalization of a theorem of Gonchar}

\begin{document}

\author{Peter Pflug}
\address{Peter Pflug\\
Carl von  Ossietzky Universit\"{a}t Oldenburg \\
Fachbereich  Mathematik\\
Postfach 2503, D--26111\\
 Oldenburg, Germany}
\email{pflug@mathematik.uni-oldenburg.de}

\author{Vi{\^e}t-Anh  Nguy\^en}
\address{Vi{\^e}t-Anh  Nguy\^en\\
Mathematics Section\\
The Abdus Salam international centre
 for theoretical physics\\
Strada costiera, 11\\
34014 Trieste, Italy}
\email{vnguyen0@ictp.trieste.it}

\subjclass[2000]{Primary 32D15, 32D10}

\keywords{Holomorphic extension, plurisubharmonic measure.}

\begin{abstract}
Let $X,\ Y$  be two complex manifolds,
  let $D\subset X,$ $ G\subset Y$ be two nonempty open sets, let
  $A$ (resp. $B$) be an open subset of  $\partial D$ (resp.
  $\partial G$), and let $W$ be the  $2$-fold cross $((D\cup A)\times B)\cup (A\times(B\cup G)).$
  Under a geometric condition on the boundary sets $A$ and $B,$
   we  show that every function  locally bounded, separately
  continuous on $W,$
continuous on $A\times B,$ and  separately
  holomorphic
  on $(A\times G) \cup (D\times B)$ ``extends" to a function  continuous on a ``domain of holomorphy" $\widehat{W}$
 and  holomorphic  on the interior of $\widehat{W}.$
\end{abstract}
\maketitle

\section{Introduction}

In the works \cite{go1,go2}  Gonchar has proved  the following remarkable  result.

\renewcommand{\thethmspec}{Gonchar's Theorem}
\begin{thmspec}
Let  $D,\ G\subset \C$ be  Jordan domains and  $A $ (resp. $B$) a nonempty open set
of the boundary $\partial D$  (resp. $\partial G$).
Let \begin{equation*}
f:\ W:=A\times (G\cup B)\bigcup (D\cup A)\times B\longrightarrow\C
\end{equation*}
 be a continuous
function  such that $f(a,\cdot)|_G$  and $f(\cdot,b)|_D$  are holomorphic for all $a\in A$ and $b\in B.$
Then there is a unique function
$\hat{f}$ continuous on
\begin{equation*}
\widehat{W} :=\left\lbrace (z,w)\in (D\cup A)\times (G\cup B) :\  \omega(z,A,D)+\omega(w,B,G)<1
\right\rbrace,
\end{equation*}
 and holomorphic on
 \begin{equation*}
 \widehat{W}^{\text{o}}:=\left\lbrace (z,w)\in D\times G :\  \omega(z,A,D)+\omega(w,B,G)<1
\right\rbrace,
 \end{equation*}
 such that $\hat{f}=f$ on $W,$  where
$\omega(\cdot,A,D),$ $\omega(\cdot,B,G)$ are the harmonic measures  (\it{See} Subsection 2.2 below).
 Moreover, if $\vert f\vert_W<\infty$ then
\begin{equation*}
 \vert \hat{f}(z,w)\vert\leq \vert f\vert_{A\times B}^{1-\omega(z,A,D)-\omega(w,B,G)} \vert
 f\vert_W^{\omega(z,A,D)+\omega(w,B,G)},\qquad (z,w)\in\widehat{W},
\end{equation*}
where $\vert g\vert_M:=\sup\limits_{M}\vert g\vert$ for a function $g$ defined on a set $M.$
\end{thmspec}

Gonchar's Theorem
generalizes   the pioneer  work of  Malgrange--Zerner \cite{ze}
on   a boundary version of the cross
theorem,
  and other results
 obtained by
Komatsu \cite{ko} and Dru\.{z}kowski \cite{dr}.
At the same time of Gonchar's work in \cite{go1},  Airapetyan and Henkin published  a  version of
 the edge-of-the-wedge theorem for CR manifolds ({\it see}  \cite{ah1}  for a brief version
 and  \cite{ah2} for a complete proof).
   Gonchar's Theorem could be deduced from the latter result.

 Recently, the authors have been able to  generalize Gonchar's result  to the case where $D,\ G$ are pseudoconvex domains
in $\C^n$  ({\it see} \cite{pn1}).

The main goal of the present  work is to
establish  a generalization of Gonchar's Theorem
for the case where $D$ and $G$ are open subsets of arbitrary complex manifolds
and $A\subset \partial D,$  $B\subset\partial G$ are open (boundary) subsets.






The proof of the result presented in this work is based on Gonchar's Theorem,
 the techniques introduced in our previous work \cite{pn1}, the
approach ``Poletsky theory of holomorphic discs and Rosay's Theorem" developed in a
recent article of the second author \cite{nv},
  and a thorough geometric study
of the plurisubharmonic measure.


\indent{\it{\bf Acknowledgment.}}   The paper was written while  the second author was visiting the
Max-Planck Institut f\"{u}r Mathematik in Bonn and the  Abdus Salam International Centre
 for Theoretical Physics
in Trieste. He wishes to express his gratitude to these organizations.

\section{Statement of the  main result and outline of the proof}

In order to state the main result, we need to introduce some notation and terminology.
 In fact, we keep the main notation from the previous work \cite{pn1}

\subsection{Topological hypersurfaces in a complex manifold}


For every open subset $U\subset \R^{2n-1}$ and every continuous function $h:\  U\longrightarrow\R,$
the graph  $\left\lbrace z=(z^{'},z_n)=(z^{'},x_n+iy_n)\in \C^n:\  y_n=h(z^{'},x_n)     \right\rbrace$
is called a {\it topological hypersurface in $\C^n$.}

 Let  $X$ be a complex manifold of dimension $n.$
 A subset $A\subset X$ is said to be  a  {\it topological hypersurface} if,
 for every point $a\in A,$ there is a    local chart  $(U,\phi:\ U\rightarrow\C^n)$ around $a$
 such that  $\phi(A\cap U)$ is a topological hypersurface in $\C^n$

Now let    $D\subset X$  be an open subset and let  $A\subset\partial D$
be an open subset
(with respect to the topology induced on $\partial D$). Suppose in addition that
$A$ is a topological hypersurface.
A point  $a\in A$ is said to be {\it   of type 1  (with respect to $D$)}
if, for every neighborhood $U$ of $a$ there is an open  neighborhood $V$ of $a$ such that $V\subset U$ and  $V\cap D$ is
a domain.  Otherwise, $a$ is said to be {\it   of type 2}.
We see easily that if $a$ is of type 2, then     for every neighborhood $U$ of $a,$
   there are an open neighborhood $V$ of $a$ and
two  domains $V_1,$ $V_2$  such that $V\subset U,$    $V\cap D=V_1\cup V_2$
 and all points in $A\cap V$ are of type 1 with respect to $V_1$ and $V_2.$

We conclude this subsection with a simple example which  may clarify  the above definitions. 
Let $G$ be the open square in $\C$
whose four vertices are $1+i,$ $-1+i,$ $-1-i,$ and $1-i.$
Define the domain
\begin{equation*}
D:=G\setminus \left [-\frac{1}{2},\frac{1}{2}\right].
\end{equation*}
Then    $A:=\partial G\cup  \left
(-\frac{1}{2},\frac{1}{2}\right)$  is not only  an open subset of $\partial D,$  but also a topological hypersurface.
 Every point of $\partial G$ is of type
1 and every point of $\left(-\frac{1}{2},\frac{1}{2}\right)$ is of type
2  (with respect to $D$).
\subsection{Plurisubharmonic measure}

Let $X$  be a complex manifold  and let  $D$ be an  open subset of $X.$
  For every function $u:\ D\longrightarrow [-\infty,\infty),$ let
\begin{equation*}
 \hat{u}(z):=
\begin{cases}
u(z),
  & z\in   D,\\
 \limsup\limits_{w\in D,\ w\to z}u(w), & z \in \partial D.
\end{cases}
\end{equation*}
For a set  $A\subset \overline{D}$ put
\begin{equation*}
h_{A,D}:=\sup\left\lbrace u\ :\  u\in\mathcal{PSH}(D),\ u\leq 1\ \text{on}\ D,\
   \hat{u}\leq 0\ \text{on}\ A    \right\rbrace,
\end{equation*}
where $\mathcal{PSH}(D)$ denotes the set   of all functions  plurisubharmonic
on $D.$

 The {\it plurisubharmonic measure} of $A$ relative to $D$ is given
 by
\begin{equation}
 \omega(z,A,D):=
   \widehat{h^{\ast}_{A,D}}(z)    , \qquad z \in  D\cup A,
\end{equation}
where $u^{\ast}$ denotes the  upper semicontinuous regularization of a function $u.$

Geometric  properties  of the plurisubharmonic measure will be discussed in   Section
3 below.

\subsection{Cross and separate holomorphicity}
Let $X,\ Y$  be two complex manifolds,
  let $D\subset X,$ $ G\subset Y$ be two nonempty open sets, let
  $A$ (resp. $B$) be either  an open subset of  $\partial D$ (resp.
  $\partial G$) or an open subset of $D$   (resp. $G$). If, moreover,  $A$  (resp.  $B$) is an open subset of
  $\partial D$  (resp.  $\partial G$),  then we assume in addition that
  $A$ (resp. $B$)  is a  topological hypersurface.

 We define
a {\it $2$-fold cross} $W,$  its {\it  interior} $W^{\text{o}},$ as
\begin{eqnarray*}
W &=&\X(A,B; D,G)
:=((D\cup A)\times B)\cup (A\times(B\cup G)),\\
W^{\text{o}} &=&\X^{\text{o}}(A,B; D,G)
:= (A\times  G)\cup (D\times B).
\end{eqnarray*}
Moreover, put
\begin{equation*}
\omega(z,w):=\omega(z,A,D)+\omega(w,B,G),\qquad
(z,w)\in (D\cup A)\times (G\cup B).
\end{equation*}

For a $2$-fold cross $W :=\X(A,B; D,G)$
define  {\it its wedge}
\begin{equation*}
\widehat{W}:=\widehat{\X}(A,B;D,G)
:=\left\lbrace (z,w)\in (D\cup A)\times (G\cup B):\ \omega(z,w)  <1
\right\rbrace.
\end{equation*}
Then the set of all interior points of the wedge
$\widehat{W}$ is given by
\begin{equation*}
 \widehat{W}^{\text{o}} :=\widehat{\X}^{\text{o}}(A,B;D,G)
 :=\left\lbrace (z,w)\in D\times G :\  \omega(z,w)<1
\right\rbrace.
\end{equation*}

We say that a function $f:W\longrightarrow \C$ is {\it separately holomorphic}
{\it on $W^{\text{o}} $} and write $f\in\mathcal{O}_s(W^{\text{o}}),$   if
 for any $a\in A $ (resp.  $b\in B$)
 the function $f(a,\cdot)|_{G}$  (resp.  $f(\cdot,b)|_{D}$ )  is holomorphic  on $G$  (resp. on $D$).

 We say that a function $f:\ W\longrightarrow \C$
   is  {\it separately continuous}
{\it on $W $}   and
write $f\in\mathcal{C}_s(W),$   if
 for any $a\in A $ (resp.  $b\in B$)
 the function $f(a,\cdot)|_{G\cup B}$  (resp.  $f(\cdot,b)|_{D\cup A}$ )  is continuous  on $G\cup B$  (resp. on $D
 \cup A$).

Throughout the paper,
 for a topological space  $M,$ $\mathcal{C}(M)$ denotes the space of  all continuous functions  $f:\ M\longrightarrow\C$
equipped with the sup-norm   $\vert f\vert_M:=\sup_M \vert f\vert.$
Moreover,  a function $f:\ M\longrightarrow\C$ is said to be {\it locally bounded} on $M$ if,
for any point $z\in M,$ there are an open neighborhood $U$ of $z$ and a
 positive number $K=K_z$ such that $\vert f\vert_U <K.$

\subsection{Statement of the main result and an outline of its proof}
We are now ready to state the  main result.

\renewcommand{\thethmspec}{Main Theorem}
  \begin{thmspec}
Let $X,\ Y$  be two complex manifolds,
  let $D\subset X,$ $ G\subset Y$ be two nonempty open sets, let
  $A$ (resp. $B$) be a nonempty open subset of  $\partial D$ (resp.
  $\partial G$).  Suppose in addition that
  $A$ and $B$ are topological hypersurfaces.
 Let  $f:\ W\longrightarrow \C$ be   such that:
\begin{itemize}
\item[ (i)]  $f\in\mathcal{C}_s(W)\cap \mathcal{O}_s(W^{\text{o}});$
\item[(ii)]   $f$ is locally bounded on $W;$
\item[ (iii)]   $f|_{A\times B}$ is continuous.
\end{itemize}

Then  there exists a unique function
$\hat{f}\in\mathcal{C}(\widehat{W})\cap \mathcal{O}(\widehat{W}^{\text{o}})$
such that  $\hat{f}=f$ on $W.$ Moreover, if $\vert f\vert_W<\infty,$ then
\begin{equation*}
 \vert \hat{f}(z,w)\vert\leq \vert f\vert_{A\times B}^{1-\omega(z,w)} \vert
 f\vert_W^{\omega(z,w)},\qquad (z,w)\in\widehat{W}.
\end{equation*}
\end{thmspec}

It is worthy to remark that the formulation of the Main Theorem bears a flavor of
 Dru\.{z}kowski's Theorem in  \cite{dr}.
In fact, when $D$ and $G$ are Jordan domains in $\C,$
 the Main Theorem follows from Gonchar's Theorem and the proof of   Dru\.{z}kowski's Theorem.
Now we give some ideas how to prove the Main Theorem.

\smallskip

In order to tackle ``arbitrary" complex manifolds,  the first  key technique here is to
apply the beautiful theorem of Rosay
\cite{ro}. This   ``Poletsky theory of holomorphic discs" approach
 has been explored in the work \cite{nv}, where  the second
author succeeded in removing the ``pseudoconvex hypothesis" in the classical cross theorems. The second key technique is to
 apply a   mixed cross type theorem ({\it see} also \cite{pn1}).
 The third key technique is to use {\it level sets} of the plurisubharmonic measure ({\it see} \cite{nv,pn1}). More precisely,
we exhaust $D$  (resp.  $G$) by the  level sets of the harmonic measure
$\omega(\cdot,A,D)$  (resp. $\omega(\cdot,B,G)$),  i.e. by
$D_{\delta}:=\left\lbrace z\in D:\ \omega(z,A,D)<1-\delta \right\rbrace$
(resp.  $G_{\delta}:=\left\lbrace w\in G:\ \omega(w,B,G)<1-\delta \right\rbrace              $)  for  $0<\delta<1.$


 Our method consists of three steps.
In the first step we  suppose that  $G$ is a  domain in $\C^m$ and   $A$ is an open subset of $D.$
In the second step we treat the case  where the pairs  $(D,A)$ and $(G,B)$ are   ``good" enough in the sense of the slicing method.
In the last one we consider the general case. For the first step we  combine the  mixed cross  theorem with
the technique of holomorphic discs.  For the second step
 one  applies the slicing method and Gonchar's Theorem.  The general philosophy is to
prove the Main Theorem with $D$  (resp. $G$) replaced by $D_{\delta}$  (resp.  $G_{\delta}$). Then we construct
the solution  for the  original open sets $D$ and $G$ by means of a gluing procedure ({\it see} also \cite{nv,pn1}).
In the last step we
transfer the holomorphicity from local situations  to  the global context.

Although our results have been stated only  for the case  of a $2$-fold cross,
they can be formulated   for the general case of an $N$-fold cross  with $N\geq 2$  ({\it see} also \cite{nv,pn1}).

\section{Preparatory results}

We present here the auxiliary  results needed for the proof of the Main Theorem.
\subsection{Poletsky theory of  discs and Rosay's Theorem on
 holomorphic  discs}
  Let $E$ denote as usual the unit disc in $\C.$ %
For a complex
manifold $\mathcal{M},$ let $\mathcal{O}(\overline{E},\mathcal{M})$ denote
the set of all holomorphic mappings $\phi:\ E\longrightarrow \mathcal{M}$ which
extend holomorphically  to   a neighborhood of  $\overline{E}.$
Such a mapping $\phi$ is called a {\it holomorphic disc} on $\mathcal{M}.$ Moreover, for
a subset $A$ of $\mathcal{M},$ let
\begin{equation*}
 1_{  A,\mathcal{M}}(z):=
\begin{cases}
1,
  &z\in   A,\\
 0, & z\in \mathcal{M}\setminus A.
\end{cases}
\end{equation*}

In the work \cite{ro}  Rosay proved the following remarkable result.
\begin{thm}\label{Rosaythm}
Let $u$ be an upper semicontinuous function on a complex manifold
$\mathcal{M}.$ Then the Poisson functional of $u$  defined by
\begin{equation*}
\mathcal{P}[u](z):=\inf\left\lbrace\frac{1}{2\pi}\int\limits_{0}^{2\pi} u(\phi(e^{i\theta}))d\theta:  \
\phi\in   \mathcal{O}(\overline{E},\mathcal{M}), \ \phi(0)=z
\right\rbrace,
\end{equation*}
is plurisubharmonic on $\mathcal{M}.$
\end{thm}

 Special cases of Theorem
\ref{Rosaythm} have been considered, for the first times, by Poletsky in \cite{po1,po2}),
and then  by L\'arusson--Sigurdsson ({\it see} \cite{ls}) and Edigarian ({\it see} \cite{ed}).

The following result  is an immediate consequence of Rosay's Theorem.
\begin{prop}\label{prop_Rosay}
Let $\mathcal{M}$ be a complex manifold and   $A$   a nonempty open
subset of $\mathcal{M}.$    Then
   $ \omega(z,A,\mathcal{M}) = \mathcal{P}[1_{\mathcal{M}\setminus A,\mathcal{M}}](z),$ $z\in\mathcal{M}.$
\end{prop}
\begin{proof}  See, for example, the proof of Proposition 3.4 in \cite{nv}.
\end{proof}

The following result is  simple but very useful.
\begin{lem}\label{lem_estimates}
Let $T$ be an open subset of $\overline{E}.$   Then
\begin{equation*}
\omega(0,T\cap E,E)\leq \frac{1}{2\pi}\int\limits_{0}^{2\pi}1_{  \partial E\setminus T,\partial E}
(e^{i\theta})d\theta.
\end{equation*}
\end{lem}
\begin{proof}  See Lemma 3.3 in \cite{nv}.
\end{proof}

\subsection{Slicing method}
 Let  $X$ be a complex manifold of dimension $n,$
 let    $D$   be an open subset of $X$  and let  $A\subset\partial D$
be an open boundary subset
which is also a topological hypersurface. We  like to study the ``slicing"  property of
$D$ near an arbitrary point  $a\in A.$  Since  our study is local, we  look at
  sufficiently small open neighborhoods $V$ of $a$  such that
$V$ is contained in a chart. Therefore, $V\cap D$  may be identified with an open neighborhood of $0$ in $\C^n.$
In addition, we may choose $V$ so that there are
 an  open subset $U\subset \R^{2n-1}$ and a
 continuous function $h:\  U\longrightarrow\R$  such that
 \begin{equation*}
 V\cap \partial D=V\cap A=\left\lbrace z=(z^{'},z_n)=(z^{'},x_n+iy_n)\in \C^n:\  y_n=h(z^{'},x_n)     \right\rbrace.
 \end{equation*}

 Assume without loss of generality that
$[-1,1]^{2n-1}\subset U $  and  $a:=0\in\C^n.$ By shrinking $U$ and $V$ (if necessary), and using the continuity of $h,$
we may find an $\epsilon>0$ such that there are only two following cases: 
\\
{\bf Case 1:} all points of $A\cap V$ are of type 1 and
\begin{equation*}
E_{z^{'}}:=\left\lbrace z=(z^{'},z_n)=(z^{'},x_n+iy_n)\in \C^n:\  -\epsilon< y_n-h(z^{'},x_n) <0    \right\rbrace\subset D\cap V
\end{equation*}
for all $  z^{'}\in  (-1,1)^{2n-2}.$\\
{\bf Case 2:} all points of $A\cap V$ are of type 2 and
\begin{multline*}
E_{z^{'}}:=\left\lbrace z=(z^{'},z_n)=(z^{'},x_n+iy_n)\in \C^n:\ \right.\\
\left. -\epsilon< y_n-h(z^{'},x_n) <0\quad
 \text{or}\quad   0< y_n-h(z^{'},x_n) <\epsilon    \right\rbrace \subset D\cap V
\end{multline*}
for all $   z^{'}\in  (-1,1)^{2n-2}.$

For a subset $S\subset\C^n,$  let $(S)_n$ denote the image of $S$ under the canonical  projection of $\C^n$
onto the $n$-th coordinate.
Observe that in  Case 1, $(E_{z^{'}})_n$ is a Jordan domain, but in  Case 2,
   $(E_{z^{'}})_n$ is a disjoint union of two Jordan domains.
  For all $z^{'}\in   (-1,1)^{2n-2},$ let
  \begin{equation*}
(A\cap V)_{z^{'}}:=\left\lbrace z=(z^{'},z_n)=(z^{'},x_n+iy_n)\in \C^n,\  :\   y_n=h(z^{'},x_n)     \right\rbrace.
\end{equation*}

\begin{defi}\label{canonical_form}
Under the above hypothesis and notation,  $(D\cap V,A\cap V)$ is said to be
  a good pair.
\end{defi}
In summary, we have shown that
\begin{prop}\label{prop_canonical_form}
 Let  $X$ be a complex manifold,
 let    $D$   be an open subset of $X$  and let  $A\subset\partial D$
be an open boundary subset
which is also a topological hypersurface. Then for all  points  $a\in A,$  there is an open neighborhood $V$ of $a$
such that  the pair $(D\cap V,A\cap V)$ is  good.
\end{prop}

  Using the construction above and the continuity of $h,$  we may apply  Rad\'o's Theorem
  (see  Theorem 2 in \cite[p. 59]{go}). Consequently, the family of harmonic measures $\omega(\cdot, ((A\cap V)_{z^{'}})_n,
  (E_{z^{'}})_n)$
  depends ``continuously" on the parameter  $z^{'}\in    (-1,1)^{2n-2}.$
 Therefore, we obtain the following result
    \begin{prop}\label{prop_slicing}
We keep the above hypothesis and notation. Consider    the set
  \begin{equation*}
  (A\cap V)_{\delta}:=\inte\Big (\bigcup\limits_{z^{'}\in    (-1,1)^{2n-2}}
 \left\lbrace z=(z^{'},z_n),\  z_n\in ( E_{z^{'}})_n:\  \omega(z_n, ((A\cap V)_{z^{'}})_n,(E_{z^{'}})_n) <\delta
   \right\rbrace \Big),
   \end{equation*}
    where  $\inte S$ denotes the set of all interior points of $S\subset \C^n.$
    Then $ (A\cap V)_{\delta}\cup (A\cap V)$ is a neighborhood of $A\cap V$ in $D\cup (A\cap V).$
     \end{prop}

   The above result has several useful consequences.

   \begin{prop}\label{prop_level}
    Let  $X$ be a complex manifold,
 let    $D$  be an open subset of $X$  and let  $A\subset\partial D$
be an open boundary subset
which is also a topological hypersurface. For all $0\leq\epsilon <1,$ let
\begin{equation*}
D_{\epsilon}:=\left\lbrace   z\in D:\  \omega(z,A,D)<1-\epsilon  \right\rbrace.
\end{equation*}
Then:\\
 1) $A$ is also  an open set of $\partial  D_{\epsilon}$  and $\lim\limits_{z\to\zeta}\omega(z, A,D)=0$ for all $\zeta\in A.$\\
 2) Moreover,
\begin{equation*}
\omega(z,A,D_{\epsilon})=\frac{\omega(z,A,D)}{1-\epsilon},\qquad  z\in D_{\epsilon}.
\end{equation*}
3)  (The Uniqueness Theorem) If $f\in\mathcal{O}(D_{\epsilon})$ be such that  $\lim\limits_{z\to\zeta} f(z)=0$ for all $\zeta\in A,$
then $f\equiv 0.$
      \end{prop}
      \begin{proof}
      Applying Proposition \ref{prop_slicing} locally,  the first assertion follows.
      Using the first assertion, the second and third ones follow from some standard arguments.
      \end{proof}

\subsection{A mixed cross theorem}
\begin{thm}\label{mixedcrossthm}
Let  $D$ be the unit disc in $\C,$ let  $G$ be  an open subset  in $\C^m,$
       let $ A$    be an open  subset of $D,$  and let $B $ be an open  subset of $\partial G.$
       Suppose in addition that $B$ is a topological hypersurface.
        Put $W:=\X(A,B;D,G)$ and
$\widehat{W}:=\widehat{\X}(A,B;D,G).$
 Let  $f:\ W\longrightarrow \C$ be   such that:
\begin{itemize}
\item[ (i)]  $f\in\mathcal{C}_s(W)\cap \mathcal{O}_s(W^{\text{o}});$
\item[(ii)]   $f$ is locally bounded on $W;$
\item[ (iii)]   $f|_{A\times B}$ is continuous.
\end{itemize}

Then  there exists a unique function
$\hat{f}\in\mathcal{C}(\widehat{W})\cap \mathcal{O}(\widehat{W}^{\text{o}})$
such that  $\hat{f}=f$ on $W.$
 Moreover, if $\vert f\vert_W<\infty,$ then
\begin{equation*}
 \vert \hat{f}(z,w)\vert\leq \vert f\vert_{A\times B}^{1-\omega(z,w)} \vert
 f\vert_W^{\omega(z,w)},\qquad (z,w)\in\widehat{W}.
\end{equation*}
\end{thm}
\begin{proof} Using Proposition \ref{prop_level},
the proof of Theorem 4.1 and Theorem 4.2 in \cite{pn1} also works in the present context
making the obviously necessary changes. In fact, the hypothesis on $D$ (i.e. $D=E$) implies that
$D$ is pseudoconvex. Therefore, we are able to apply    the classical  method of doubly orthogonal bases of Bergman type.
\end{proof}

\section{Part 1 of the proof of the Main Theorem}

The main purpose of this section is  to prove the following
mixed cross theorem.
\begin{thm}\label{thm_Part1}
Let  $D$ be a  complex manifold, let    $G$ be  an open subset  in $\C^m,$
       let $ A$    be an open  subset of $D,$  and let $B $ be an open  subset of $\partial G.$
       Suppose in addition that $B$ is a topological hypersurface.
        Put $W:=\X(A,B;D,G)$ and
$\widehat{W}:=\widehat{\X}(A,B;D,G).$
 Let  $f:\ W\longrightarrow \C$ be   such that:
\begin{itemize}
\item[ (i)]  $f\in\mathcal{C}_s(W)\cap \mathcal{O}_s(W^{\text{o}});$
\item[(ii)]   $f$ is locally bounded on $W;$
\item[ (iii)]   $f|_{A\times B}$ is continuous.
\end{itemize}

Then  there exists a unique function
$\hat{f}\in\mathcal{C}(\widehat{W})\cap \mathcal{O}(\widehat{W}^{\text{o}})$
such that  $\hat{f}=f$ on $W.$
 Moreover, if $\vert f\vert_W<\infty,$ then
\begin{equation*}
 \vert \hat{f}(z,w)\vert\leq \vert f\vert_{A\times B}^{1-\omega(z,w)} \vert
 f\vert_W^{\omega(z,w)},\qquad (z,w)\in\widehat{W}.
\end{equation*}
\end{thm}
It is worthy to remark that  Theorem \ref{thm_Part1} removes the hypothesis ``pseudoconvex" in Theorem  \ref{mixedcrossthm}.
\begin{proof}
It follows essentially the proof of Theorem 4.1 in \cite{nv}.   We begin the proof with the following lemma.
\begin{lem}
\label{lem_Part1}
We keep the hypothesis of Theorem \ref{thm_Part1}. For $j\in \{1,2\},$ let $\phi_j\in\mathcal{O}(
 \overline{E}, D) $ be a holomorphic disc,  and let $t_j\in E$ such that $\phi_1(t_1)=\phi_2(t_2)$ and
$  \frac{1}{2\pi}\int\limits_{0}^{2\pi}1_{D\setminus A,D} (\phi_j(e^{i\theta}))d\theta
<1,$ $j=1,2.$ Then:
\\
1) For  $j \in \{
1,2\},$  the function $(t,w)\mapsto f(\phi_j(t),w)$ belongs to
\begin{equation*}
 \mathcal{C}_s(\X(\phi^{-1}_j(A)\cap  E,B;E,G))\cap
 \mathcal{O}_s(\X^{\text{o}}(\phi^{-1}_j(A)\cap  E,B;E,G)),
 \end{equation*}
 and continuous on  $(\phi^{-1}_j(A)\cap  E)\times  B,$ where
$\phi^{-1}_j(A):=\lbrace t\in \overline{E}:\ \phi_j(t)\in A\rbrace.$     \\
2)  For  $j \in \{
1,2\},$ in virtue of Part 1) and applying  Theorem   \ref{mixedcrossthm},
 let $\hat{f}_j$ be the unique function in
\begin{equation*}
 \mathcal{C}\left(\widehat{\X}(\phi^{-1}_j(A)\cap  E,B;E,G)\right)
 \cap  \mathcal{O}\left(\widehat{\X}^{\text{o}}(\phi^{-1}_j(A)\cap  E,B;E,G)\right)
 \end{equation*}
  such that  $\hat{f}_j(t,w)=f(\phi_j(t),w),$
 $(t,w)\in\X\left(\phi^{-1}_j(A)\cap  E,B;E,G\right).$   Then
\begin{equation*}
\hat{f}_1(t_1,w)=\hat{f}_2(t_2,w),
\end{equation*}
for all $w\in G$ such that
$(t_j,w)\in\widehat{\X}\left(\phi^{-1}_j(A)\cap  E,B;E,G\right),$  $j \in \{1,2\}.$
\end{lem}

\smallskip
\noindent {\it Proof of Lemma \ref{lem_Part1}.}
 Part 1) follows immediately from the hypothesis. Therefore, it remains to prove Part
 2). To do this fix  $w_0\in G\cup B$ such that $(t_j,w_0)\in\widehat{\X}\left(\phi^{-1}_j(A)\cap E
 ,B;E,G\right)$
 for  $j \in \{1,2\}.$
We need to show that $\hat{f}_1(t_1,w_0)=\hat{f}_2(t_2,w_0).$  Observe
that both functions  $w\in\mathcal{G}\mapsto \hat{f}_1(t_1,w)$ and $w\in\mathcal{G}\mapsto\hat{f}_2(t_2,w)$ belong to  $\mathcal{O}
(\mathcal{G}),$
where $\mathcal{G}$ is the connected component which contains $w_0$ of the
following open set
\begin{equation*}
  \left\lbrace w\in G:\ \omega(w,B,G)<1-\max\limits_{ j \in \{1,2\}}  \omega(t_j,\phi^{-1}_j(A)\cap E,E)
   \right\rbrace.
   \end{equation*}
    On the other hand,    for any  $j \in \{
1,2\}$  and $w\in B,$
   $(t_j,w)\in\widehat{\X}\left(\phi^{-1}_j(A)\cap  E,B;E,G\right).$   This, combined with the equality
   $\phi_1(t_1)=\phi_2(t_2),$ implies that
\begin{equation*}
\hat{f}_1(t_1,w)=f(\phi_1(t_1),w)=f(\phi_2(t_2),w)=\hat{f}_2(t_2,w) , \qquad w\in B.
\end{equation*}
Therefore, by the Uniqueness Theorem ({\it  see} Part 3) of Proposition \ref{prop_level}),   $\hat{f}_1(t_1,w)=\hat{f}_2(t_2,w),$
  $  w\in\mathcal{G}.$ Hence,  $\hat{f}_1(t_1,w_0)=\hat{f}_2(t_2,w_0),$
which completes the proof of the lemma.
\hfill $\square$

\medskip

Now we return to the proof of the theorem.
 We define $\hat{f}$ as follows: Let $\mathcal{W}$ be  the set of all pairs
$(z,w)\in D\times (G\cup B)$  with the property that there are a holomorphic disc
$\phi\in\mathcal{O}(\overline{E},D)$ and $t\in E$ such that $\phi(t)=z$
and  $(t,w)\in\widehat{\X}\left(\phi^{-1}(A)\cap  E,B;E,G\right).$  In virtue of
  Theorem  \ref{mixedcrossthm},
 let $\hat{f}_{\phi}$ be the unique function  in
\begin{equation*}
 \mathcal{C}\left(\widehat{\X}(\phi^{-1}(A)\cap  E,B;E,G)\right)\cap \mathcal{O}\left(\widehat{\X}^{\text{o}}
 (\phi^{-1}(A)\cap  E,B;E,G)\right)
 \end{equation*}
  such that
\begin{equation}\label{eq4.1}
  \hat{f}_{\phi}(t,w)=f(\phi(t),w),\qquad (t,w)\in \X\left(\phi^{-1}(A)\cap  E,B;E,G\right).
\end{equation}
   Then  the desired extension function $\hat{f}$ is given by
\begin{equation}\label{eq4.2}
 \hat{f}(z,w):=\hat{f}_{\phi}(t,w) .
\end{equation}
In virtue of Part 2) of Lemma \ref{lem_Part1}, $\hat{f}$ is well-defined on $\mathcal{W}.$
We next prove that
\begin{equation}
 \label{eq4.3}
 \mathcal{W}=\widehat{W}  .
\end{equation}
Taking  (\ref{eq4.3}) for granted, then $\hat{f}$ is
well-defined on $\widehat{W}.$  

Now we return to (\ref{eq4.3}). To prove the inclusion $\mathcal{W}\subset\widehat{W},$
let  $(z,w)\in\mathcal{W}.$  By the above definition of $\mathcal{W},$ one may find a  holomorphic  disc    $\phi\in
\mathcal{O}(\overline{E}, D),$  a point  $t\in E$ such that $\phi(t)=z$
and  $(t,w)\in\widehat{\X}\left(\phi^{-1}(A)\cap  E,B;E,G\right).$
Since  $ \omega(\phi(t),A,D)\leq  \omega(t,\phi^{-1}(A)\cap  E,E),$ it follows that
\begin{equation*}
 \omega (z,A,D) +  \omega(w,B,G) \leq   \omega(t,\phi^{-1}(A)\cap  E,E)+
  \omega(w,B,G)<1,
\end{equation*}
Hence  $(z,w)\in\widehat{W}.$ This proves the above mentioned inclusion.

To finish the proof of  (\ref{eq4.3}), it suffices to show that  $\widehat{W}\subset\mathcal{W}.$
To do this, let $(z,w)\in\widehat{W}$ and fix any $\epsilon>0$ such that
\begin{equation}\label{eq4.4}
\epsilon<1-  \omega(z,A,D) -  \omega(w,B,G).
\end{equation}
 Applying  Theorem \ref{Rosaythm} and Proposition \ref{prop_Rosay}, there is a holomorphic  disc   $\phi\in
\mathcal{O}(\overline{E},  D)$   such that $\phi(0)=z$ and
\begin{equation}\label{eq4.5}
\frac{1}{2\pi} \int\limits_{0}^{2\pi} 1_{D\setminus
A,D}(\phi(e^{i\theta}))d\theta< \omega(z,A,D)+\epsilon.
\end{equation}
Observe that
\begin{eqnarray*}
  \omega(0,\phi^{-1}(A)\cap  E,E)+ \omega(w,B,G)&\leq &\frac{1}{2\pi} \int\limits_{0}^{2\pi} 1_{D\setminus
A,D}(\phi(e^{i\theta}))d\theta+ \omega(w,B,G)\\
 &< &  \omega(z,A,D) +  \omega(w,B,G)+\epsilon<1,
\end{eqnarray*}
where the first inequality follows from Lemma \ref{lem_estimates},
the second one   from (\ref{eq4.5}), and the last one from (\ref{eq4.4}). Hence,
 $(0,w)\in\widehat{\X}\left(\phi^{-1}(A)\cap  E,B;E,G\right),$ which
 implies that $(z,w)\in\mathcal{W}.$ This complete  the proof of  (\ref{eq4.3}).
Hence, the construction of the extension function $\hat{f}$ on $\widehat{W}$ has been completed.

Using  (\ref{eq4.1})--(\ref{eq4.3}), the proof given in Step 2 and 3 of Section 4 in \cite{nv} still works in the present context making the
obviously necessary changes. This gives that $\hat{f}=f$ on $W$ and $\hat{f}\in \mathcal{C}(\widehat{W})\cap\mathcal{O}(\widehat{W}^{\text{o}}).$
Consequently,     arguing as in the proof of Theorem 4.2 in \cite{pn1}, the desired estimate of the theorem follows.


This completes the proof of Theorem  4.1.
\end{proof}
\section{Part 2 of the proof: Local result}

The main purpose of this section is to prove  the following ``local" result.
\begin{thm}\label{thm_Part2}
Let $(D,A)$ and $(G,B)$ be two good pairs.
 Let  $f:\ W\longrightarrow \C$ be   such that:
\begin{itemize}
\item[ (i)]  $f\in\mathcal{C}_s(W)\cap \mathcal{O}_s(W^{\text{o}});$
\item[(ii)]   $f$ is locally bounded on $W;$
\item[ (iii)]   $f|_{A\times B}$ is continuous.
\end{itemize}

Then  there exists a unique function
$\hat{f}\in\mathcal{C}(\widehat{W})\cap \mathcal{O}(\widehat{W}^{\text{o}})$
such that  $\hat{f}=f$ on $W.$ Moreover, if $\vert f\vert_W<\infty,$ then
\begin{equation*}
 \vert \hat{f}(z,w)\vert\leq \vert f\vert_{A\times B}^{1-\omega(z,w)} \vert
 f\vert_W^{\omega(z,w)},\qquad (z,w)\in\widehat{W}.
\end{equation*}
\end{thm}
\begin{proof}
We assume without loss of generality that $D\subset \C^n,$  $G\subset \C^m.$

For every $0<\delta<\frac{1}{2},$ define
\begin{equation}\label{eq1_Part2}
\begin{split}
E_{z^{'},\delta}&:=\left\lbrace z=(z^{'},z_n):\  z_n\in (E_{z^{'}})_n,\ \omega(z_n,(A_{z^{'}})_n,(E_{z^{'}})_n)<\delta   \right\rbrace,\qquad z^{'}\in (-1,1)^{2n-2},\\
E_{w^{'},\delta}&:=\left\lbrace w=(w^{'},w_m):\  w_m\in (E_{w^{'}})_m,\ \omega(w_m,(B_{w^{'}})_m,(E_{w^{'}})_m)<\delta   \right\rbrace,\qquad w^{'}\in (-1,1)^{2m-2},\\
D_{\delta}&:=\left\lbrace z\in D:\  \omega(z,A,D)<1-\delta    \right\rbrace,\qquad
G_{\delta}:=\left\lbrace w\in G:\  \omega(w,B,G)<1-\delta    \right\rbrace,\\
A_{\delta}& :=\inte\Big(\bigcup\limits_{z^{'}\in  (-1,1)^{2n-2} } E_{z^{'},\delta}\Big),\qquad
B_{\delta} :=\inte\Big(\bigcup\limits_{w^{'}\in  (-1,1)^{2m-2}} E_{w^{'},\delta}\Big).
\end{split}
\end{equation}

The proof is divided  into two steps.\\
{\bf  Step 1:  $G$ is a Jordan domain.}

Firstly, we apply the slicing method: For all $z^{'}\in  (-1,1)^{2n-2} ,$  consider the function
\begin{equation}\label{eq1_Step1_Part2}
f_{z^{'}}(z_n,w):=f(z,w),\qquad  (z_n,w)\in \X((A_{z^{'}}\cap \partial E_{z^{'}})_n,B;(E_{z^{'}})_n,G).
\end{equation}
Applying Gonchar's Theorem, we obtain an extension function
\begin{equation*}
\hat{f}_{z^{'}}\in\mathcal{C}\Big( \widehat{\X}((A_{z^{'}}\cap \partial E_{z^{'}})_n,B;(E_{z^{'}})_n,G)
  \Big)\cap \mathcal{O}\Big(
\widehat{ \X}^{\text{o}}((A_{z^{'}}\cap \partial E_{z^{'}})_n,B;(E_{z^{'}})_n,G)  \Big))
\end{equation*}
such that
\begin{equation}\label{eq2_Step1_Part2}
\hat{f}_{z^{'}}(z_n,w)=f_{z^{'}}(z_n,w),\qquad  (z_n,w)\in \X((A_{z^{'}}\cap \partial E_{z^{'}})_n,B;(E_{z^{'}})_n,G).
\end{equation}
Using (\ref{eq1_Part2})--(\ref{eq2_Step1_Part2}),
 we are able to define a new function  $\tilde{f}_{\delta}$ on $ \X(A_{\delta},B;D, G_{\delta})$  as follows
\begin{equation}\label{eq3_Step1_Part2}
 \tilde{f}_{\delta}(z,w):=
\begin{cases}
 \hat{f}_{z^{'}}(z_n,w),
  & \qquad (z,w)\in A_{\delta}\times G_{\delta}, \\
  f(z,w), &   \qquad  (z,w)\in  D\times B        .
\end{cases}
\end{equation}
Applying Theorem \ref{thm_Part1}, we obtain  an extension function
\begin{equation*}
\hat{f}_{\delta}\in\mathcal{C}\Big(\widehat{\X}(A_{\delta},B;D, G_{\delta})   \Big)\cap
 \mathcal{O}\Big(\widehat{\X}^{\text{o}}(A_{\delta},B;D, G_{\delta})   \Big)
\end{equation*}
such that
\begin{equation}\label{eq4_Step1_Part2}
\hat{f}_{\delta}(z,w)=\tilde{f}_{\delta}(z,w),\qquad  (z,w)\in \X(A_{\delta},B;D,G_{\delta}).
\end{equation}
On the other hand, using Proposition \ref{prop_level}  and (\ref{eq1_Part2}), we see that
\begin{equation} \label{eq5_Step1_Part2}
\lim\limits_{\delta\to 0^{+}}\omega(z,A_{\delta},D)=\omega(z,A,D)\qquad\text{and}\qquad
\lim\limits_{\delta\to 0^{+}}\omega(w,B,G_{\delta})=\omega(w,B,G).
\end{equation}
We are now in a position to define the desired extension  function $\hat{f}.$
Indeed, one  glues
$\left(\hat{f}_{\delta}\right)_{0<\delta<\frac{1}{2}}$ together to obtain
$\hat{f}$ in the following way
\begin{equation}\label{eq6_Step1_Part2}
\hat{f}:=\begin{cases}
\lim\limits_{\delta\to 0} \hat{f}_{\delta} &\qquad \text{on}\
 \widehat{W}^{\text{o}}\cup (D\times B),\\
   f &\qquad \text{on}\ A\times G.
   \end{cases}
\end{equation}
Using (\ref{eq1_Step1_Part2})--(\ref{eq5_Step1_Part2}) and a gluing argument as in Lemma 6.5 in \cite{nv},
 it can be checked that the limit (\ref{eq6_Step1_Part2}) exists and possesses all the required
properties.

\noindent {\bf  Step 2: The general case.}

Firstly, we apply the slicing method: For all $z^{'}\in   (-1,1)^{2n-2}$  and $w^{'}\in   (-1,1)^{2m-2},$ consider the functions
\begin{equation}\label{eq1_Step2_Part2}
\begin{split}
f_{z^{'}}(z_n,w)&:=f(z,w),\qquad  (z_n,w)\in \X((A_{z^{'}}\cap \partial E_{z^{'}})_n,B;(E_{z^{'}})_n,G),\\
f_{w^{'}}(z,w_m)&:=f(z,w),\qquad  (z,w_m)\in \X(A,(B_{w^{'}}\cap \partial E_{w^{'}})_m;D,(E_{w^{'}})_m).
\end{split}
\end{equation}
Applying the result of Step 1, we obtain extension functions
\begin{eqnarray*}
\hat{f}_{z^{'}}&\in&\mathcal{C}\Big( \widehat{\X}((A_{z^{'}}\cap \partial E_{z^{'}})_n,B;(E_{z^{'}})_n,G)
  \Big)\cap \mathcal{O}\Big(
\widehat{ \X}^{\text{o}}((A_{z^{'}}\cap \partial E_{z^{'}})_n,B;(E_{z^{'}})_n,G)  \Big),\\
\hat{f}_{w^{'}}&\in&\mathcal{C}\Big( \widehat{\X}(A,(B_{w^{'}}\cap \partial E_{w^{'}})_m;D,(E_{w^{'}})_m)
  \Big)\cap \mathcal{O}\Big(
\widehat{ \X}^{\text{o}}(A,(B_{w^{'}}\cap \partial E_{w^{'}})_m;D,(E_{w^{'}})_m)  \Big)
\end{eqnarray*}
such that
\begin{equation}\label{eq2_Step2_Part2}
\begin{split}
\hat{f}_{z^{'}}(z_n,w)&=f_{z^{'}}(z_n,w),\qquad  (z_n,w)\in \X((A_{z^{'}}\cap \partial E_{z^{'}})_n,B;(E_{z^{'}})_n,G),\\
\hat{f}_{w^{'}}(z,w_m)&=f_{w^{'}}(z,w_m),\qquad  (z,w_m)\in \X(A,(B_{w^{'}}\cap \partial E_{w^{'}})_m;D,(E_{w^{'}})_m).
\end{split}
\end{equation}
Using  (\ref{eq1_Part2})--(\ref{eq2_Step1_Part2}) and (\ref{eq1_Step2_Part2})-- (\ref{eq2_Step2_Part2}), it can be checked that
\begin{equation*}
\hat{f}_{z^{'}}(z_n,w)=f_{w^{'}}(z,w_m),\qquad  (z,w)\in A_{\delta}\times B_{\delta}.
\end{equation*}
 Therefore,
 we are able to define a new function  $\tilde{f}_{\delta}$ on $ \X(A_{\delta},B_{\delta};D_{\delta}, G_{\delta})$  as follows
\begin{equation}\label{eq3_Step2_Part2}
 \tilde{f}_{\delta}(z,w):=
\begin{cases}
 \hat{f}_{z^{'}}(z_n,w),
  & \qquad (z,w)\in A_{\delta}\times G_{\delta}, \\
  \hat{f}_{w^{'}}(z,w_n), &   \qquad  (z,w)\in  D_{\delta}\times B_{\delta}        .
\end{cases}
\end{equation}

Applying Theorem A  or Theorem 5.1 in \cite{nv}, we obtain  an extension function
$\hat{f}_{\delta}\in \mathcal{O}\Big(\widehat{\X}^{\text{o}}(A_{\delta},B_{\delta};D_{\delta}, G_{\delta})   \Big)$
such that
\begin{equation}\label{eq4_Step2_Part2}
\hat{f}_{\delta}(z,w)=\tilde{f}_{\delta}(z,w),\qquad  (z,w)\in \X(A_{\delta},B_{\delta};D_{\delta},G_{\delta}).
\end{equation}

We are now in a position to define the desired extension  function $\hat{f}.$
Indeed, one  glues
$\left(\hat{f}_{\delta}\right)_{0<\delta<\frac{1}{2}}$ together to obtain
$\hat{f}$ in the following way
\begin{equation}\label{eq5_Step2_Part2}
\hat{f}:=\begin{cases}
\lim\limits_{\delta\to 0} \hat{f}_{\delta} &\qquad \text{on}\
 \widehat{W}^{\text{o}},\\
   f &\qquad \text{on}\ W.
   \end{cases}
\end{equation}
Using  the first identity in (\ref{eq5_Step1_Part2}) and   (\ref{eq1_Step2_Part2})--(\ref{eq4_Step2_Part2})
and applying  Lemma 6.5 in \cite{nv}, we can prove that the function given by the limit (\ref{eq5_Step2_Part2})
 exists. Moreover,  $\hat{f}=f$ on $W$ and $\hat{f}\in \mathcal{C}(\widehat{W})\cap\mathcal{O}(\widehat{W}^{\text{o}}).$
Consequently,     arguing as in the proof of Theorem 4.2 in \cite{pn1}, the desired estimate of the theorem follows.
Hence, the proof is complete.
\end{proof}
\section{Proof of the Main Theorem}
%
%
%
%


By Proposition  \ref{prop_canonical_form}, for all $a$ (resp. $b\in B$) we may  fix an
open neighborhood $U_{a}$ of  $a$  (resp.  $V_b$ of $b$) such that
 $(D\cap U_{a},A\cap U_{a}) $
(resp. $(G\cap V_b,B\cap V_b$) ) is a good pair.
 For any $0<\delta<\frac{1}{2},$  define
\begin{equation}\label{eq1_Part3}
\begin{split}
U_{a,\delta}&:=\left\lbrace z\in U_{a}\cap D:\ \omega(z, A\cap U_a,   U_a\cap D)<\delta  \right\rbrace,\qquad
a\in A,\\
V_{b,\delta}&:=\left\lbrace w\in V_{b}\cap G:\  \omega(w, B\cap V_b,   V_b\cap G)<\delta  \right\rbrace,\qquad
b\in B,\\
A_{\delta}&:=\bigcup\limits_{a\in A} U_{a,\delta},\qquad
B_{\delta}:=\bigcup\limits_{b\in B} V_{b,\delta},\\
 D_{\delta}&:=\left\lbrace z\in D:\
  \omega(z,A,D)<1-\delta\right\rbrace,\quad
  G_{\delta}:=\left\lbrace w\in G:\  \omega(w,B,G)<1-\delta\right\rbrace.
\end{split}
\end{equation}

We divide the proof into two steps.\\
{\bf  Step 1: $(G,B)$ is a good pair.}

Suppose without loss of generality that  $G\subset\C^m.$
 For each $a\in A,$  let $f_a:=f|_{\X\left( A\cap U_{a} ,B;  D\cap U_{a},G\right)}.$
Using the hypothesis on $f$ we deduce that   $f_a$ is locally bounded,
\begin{equation*}
f_a\in\mathcal{C}_s\Big(\X\left( A\cap U_{a} ,B; D\cap U_{a},G\right)\Big)\cap\mathcal{O}_s\Big(\X\left( A\cap U_{a} ,B; D\cap U_{a},G\right)\Big)
\end{equation*}
and  that  $f_a|_{(A\cap U_a)\times B}\in \mathcal{C}\Big((A\cap U_a)\times B\Big).$ Recall that
$(D\cap U_a,A\cap U_a)$  and $(G,B)$ are good pairs.  Consequently,  applying Theorem
\ref{thm_Part2} to $f_a$ yields that there is a unique function
\begin{equation*}
\hat{f}_{a} \in   \mathcal{C}\Big(\widehat{\X}
 \left(A\cap U_{a} ,B;D\cap  U_{a},G\right)\Big)\cap  \mathcal{O}\Big(\widehat{\X}^{\text{o}}
 \left(A\cap U_{a} ,B;D\cap  U_{a},G\right)\Big)
\end{equation*}
such that
\begin{equation}\label{eq1_Step1_Part3}
\hat{f}_{a}(z,w)=f_a(z,w)=f(z,w),\  (z,w)\in
 \X\left(  A\cap   U_a,B;  D\cap U_{a},G\right) .
\end{equation}
Arguing as in Lemma 6.4 in \cite{nv} and using Definition 6.3 therein, we can show that
 the family $\left(\hat{f}_{a}|_{U_{a,\delta}\times G_{\delta}} \right)_{a\in A}$  is collective
 for  all  $0<\delta< \frac{1}{2}.$ In virtue of (\ref{eq1_Part3}),
let
\begin{equation}\label{eq2_Step1_Part3}
\tilde{\tilde{f}}_{\delta}\in \mathcal{O}(A_{\delta}\times
G_{\delta})
\end{equation}
denote the collected function of this family.
 In virtue of  (\ref{eq1_Step1_Part3})--(\ref{eq2_Step1_Part3}),
 we are able to define a new function $\tilde{f}_{\delta}$ on $\X\left(A_{\delta},B;D,G_{\delta}
 \right)$ as follows
\begin{equation*}
 \tilde{f}_{\delta}:=
\begin{cases}
 \tilde{\tilde{f}}_{\delta},
  & \qquad\text{on}\  A_{\delta}\times G_{\delta}, \\
  f, &   \qquad\text{on}\ D\times B        .
\end{cases}
\end{equation*}
Using this  and  (\ref{eq1_Step1_Part3})--(\ref{eq2_Step1_Part3}), we  see that  $ \tilde{f}_{\delta}\in \mathcal{O}_s\Big(
\X\left(A_{\delta}, B;D,G_{\delta}
 \right)\Big),$ and
\begin{equation}\label{eq3_Step1_Part3}
  \tilde{f}_{\delta}=f\qquad\text{on}\ D\times B.
\end{equation}
 Since $A_{\delta}$ is open in $D,$ and $B$ is not only an open set of  $\partial G_{\delta},$  but also a topological hypersurface
 (by Proposition   \ref{prop_level}),  we are able to  apply  Theorem  \ref{thm_Part1}
to  $ \tilde{f}_{\delta}$ in order to obtain a function
 \begin{equation*}
 \hat{f}_{\delta}\in \mathcal{C}\Big(  \widehat{\X}\left(A_{\delta}, B;D,G_{\delta}
 \right)\Big)\cap \mathcal{O}\Big(  \widehat{\X}^{\text{o}}\left(A_{\delta}, B;D,G_{\delta}
 \right)\Big)
 \end{equation*}  such that
\begin{equation}\label{eq4_Step1_Part3}
  \hat{f}_{\delta}=  \tilde{f}_{\delta}\qquad\text{on}\ \X\left(A_{\delta},
  B;D,G_{\delta}\right).
\end{equation}

We are now in a position to define the desired extension  function $\hat{f}.$
Indeed, one  glues
$\left(\hat{f}_{\delta}\right)_{0<\delta<\frac{1}{2}}$ together to obtain
$\hat{f}$ in the following way
\begin{equation}\label{eq5_Step1_Part3}
 \hat{f}:=\begin{cases}
\lim\limits_{\delta\to 0} \hat{f}_{\delta} &\qquad \text{on}\
 \widehat{W}^{\text{o}},\\
   f &\qquad \text{on}\ W.
   \end{cases}
 \end{equation}
Using   (\ref{eq1_Step1_Part3})--(\ref{eq5_Step1_Part3})  and arguing as in (6.12)--(6.14) in \cite{nv},
we  see that $\hat{f}$  is well-defined   and possesses all the required
properties.
\\
{\bf  Step 2:  The general case.}

 For each $a\in A,$  let $f_a:=f|_{\X\left( A\cap U_{a} ,B; D\cap U_{a},G\right)}.$
Using the hypothesis on $f$ and the fact that  $(D\cap U_a,A\cap U_a)$ is
a good pair, we are able to apply  the result of Step 1
 to $f_a.$  Consequently,  there is a unique function
 \begin{equation*}
   \hat{f}_{a} \in   \mathcal{C}\Big(\widehat{\X}
 \left(A\cap U_{a} ,B; D\cap U_{a},G\right)\Big)\cap  \mathcal{O}\Big(\widehat{\X}^{\text{o}}
 \left(A\cap U_{a} ,B; D\cap U_{a},G\right)\Big)
 \end{equation*}
such that
\begin{equation}\label{eq1_Step2_Part3}
\hat{f}_{a}(z,w)=f(z,w),\qquad  (z,w)\in
 \X\left(  A\cap  U_a,B; D\cap  U_{a},G\right) .
\end{equation}
 Let $0<\delta<\frac{1}{2}.$ In virtue of (\ref{eq1_Part3}) and  (\ref{eq1_Step2_Part3}), we may apply Lemma 6.4 in \cite{nv}.   Consequently, we can collect  the family
 $\left(\hat{f}_{a}|_{U_{a,\delta}\times G_{\delta}} \right)_{a\in A}$
in order to obtain the collected   function $\tilde{f}^A_{\delta}\in \mathcal{O}(A_{\delta}\times
G_{\delta}).$

Similarly,  for each $b\in B,$   one obtains a unique function
 \begin{equation*}
 \hat{f}_{b} \in  \mathcal{C}\Big(\widehat{\X}
 \left(A, B\cap V_{b} ; D,G\cap V_{b}\right)\Big)\cap  \mathcal{O}\Big(\widehat{\X}^{\text{o}}
 \left(A, B\cap V_{b} ; D,G\cap V_{b}\right)\Big)
 \end{equation*}
such that
\begin{equation}\label{eq2_Step2_Part3}
\hat{f}_{b}(z,w) =f(z,w),\qquad (z,w)\in
 \X\left(  A,B\cap V_b;  D,V_{b}\right) .
\end{equation}
Moreover, one can collect  the family
 $\left(\hat{f}_{b}|_{D_{\delta}\times V_{b,\delta}} \right)_{b\in B}$
in order to obtain the collected   function $\tilde{f}^B_{\delta}\in \mathcal{O}(D_{\delta}\times
B_{\delta}).$

Arguing as in the proof of (6.17)--(6.18) in \cite{nv}, we can show  that
\begin{equation*}
 \tilde{f}^A_{\delta}=\tilde{f}^B_{\delta}\qquad \text{on}\  A_{\delta}\times
B_{\delta}.
\end{equation*}
Consequently,  we are able to define a new  function $\tilde{f}_{\delta}$ on $
 \X\left(A_{\delta}, B_{\delta};D_{\delta},
 G_{\delta}\right)$ as follows
\begin{equation}\label{eq3_Step2_Part3}
 \tilde{f}_{\delta}:=
\begin{cases}
\tilde{f}^A_{\delta},
  & \qquad\text{on}\  A_{\delta}\times G_{\delta}, \\
  \tilde{f}^B_{\delta}, &   \qquad\text{on}\ D_{\delta}\times B_{\delta}        .
\end{cases}
\end{equation}
Using formula (\ref{eq3_Step2_Part3}) it can be readily checked that
$\tilde{f}_{\delta}\in \mathcal{O}_s\Big(\X\left(A_{\delta}, B_{\delta};D_{\delta},G_{\delta}
 \right)\Big).$ Since  we know that  $A_{\delta}$
 (resp.  $B_{\delta}$) is an open subset of $D_{\delta}$  (resp.
 $G_{\delta}$), we are able to apply  Theorem  A or Theorem 5.1 in \cite{nv}
to  $ \tilde{f}_{\delta}$  for every $0<\delta<\frac{1}{2}.$ Consequently,
one
  obtains a unique  function $\hat{f}_{\delta}\in \mathcal{O}\Big(
 \widehat{\X}\left(A_{\delta}, B_{\delta};D_{\delta},G_{\delta}
 \right)\Big)$ such that
\begin{equation}\label{eq4_Step2_Part3}
  \hat{f}_{\delta}=  \tilde{f}_{\delta}\qquad\text{on}\ \X\left(A_{\delta},
  B_{\delta};D_{\delta},G_{\delta}\right).
\end{equation}

We are now in a position to define the desired extension  function $\hat{f}.$
 \begin{equation}\label{eq5_Step2_Part3}
\hat{f}:=
\begin{cases}
\lim\limits_{\delta\to 0} \hat{f}_{\delta},& \qquad \text{on}\
 \widehat{W}^{\text{o}},\\
   f,&\qquad \text{on}\  W.
\end{cases}
\end{equation}
 To prove that $\hat{f}$  is well-defined and $\hat{f}=f$ on $W$ and $\hat{f}\in \mathcal{C}(\widehat{W})\cap\mathcal{O}(\widehat{W}^{\text{o}}),$
 one proceeds as in the end of the proof of Theorem 6.1
 in \cite{nv} using (\ref{eq1_Step2_Part3})--(\ref{eq5_Step2_Part3}). Consequently,     arguing as in the proof of Theorem 4.2 in \cite{pn1}, the desired estimate of the theorem follows.
Hence, the proof of the Main Theorem is complete.

\end{document}